\documentclass[11pt]{article}
\usepackage{amsmath,fullpage,amssymb,amsthm,hyperref,enumerate,shuffle}
\usepackage{tikz,color}

\theoremstyle{remark}

\theoremstyle{definition}

\newcommand{\D}{\mathcal{D}}
\newcommand{\F}{\mathcal{F}}
\renewcommand{\S}{\mathcal{S}}
\newcommand\red{\textcolor{red}}
\newcommand\blue{\textcolor{blue}}

\title{A simple bijective proof of a familiar derangement recurrence}

\author{Sergi Elizalde\thanks{Department of Mathematics, Dartmouth College, Hanover, NH 03755. \texttt{sergi.elizalde@dartmouth.edu}}}

\date{}

\begin{document}

\maketitle

\begin{abstract}
It is well known that the derangement numbers $d_n$, which count permutations of length $n$ with no fixed points, satisfy the recurrence $d_n=nd_{n-1}+(-1)^n$ for $n\ge1$. Combinatorial proofs of this formula have been given by Remmel, Wilf, D\'esarm\'enien and Benjamin--Ornstein. Here we present yet another, arguably simpler, bijective proof.
\end{abstract}

Let $\S_n$ denote the set of permutations of $\{1,2,\dots,n\}$. A fixed point of $\pi\in\S_n$ is an element~$i$ such that $\pi(i)=i$. 
Let $\D_n\subseteq\S_n$ denote the set of permutations with no fixed points, often called {\em derangements}, and let $d_n=|\D_n|$.
Let $\F_n\subseteq\S_n$ denote the set of permutations with exactly one fixed point. Clearly, $|\F_n|=nd_{n-1}$, since permutations in $\F_n$ are determined by choosing the fixed point among $n$ possibilities, and then taking a derangement of the remaining $n-1$ elements.

It is well known \cite[Eq.~(2.13)]{Stanley} that the derangement numbers $d_n$ satisfy the recurrence 
\begin{equation}\label{eq:recdn}
d_n=nd_{n-1}+(-1)^n
\end{equation}
for $n\ge1$. This equation states that the number of $\pi\in\S_n$ with no fixed points and the number of $\pi\in\S_n$ with one fixed point differ by one. Stanley~\cite{Stanley} acknowledges that proving recurrence~\eqref{eq:recdn} combinatorially requires considerably more work than proving the other well-known recurrence for derangement numbers, $d_n=(n-1)(d_{n-1}+d_{n-2})$.
%The first derangement numbers $d_n$ for $n\ge0$ are $1, 0, 1, 2, 9, 44, 265, 1854,\dots$.
Bijective proofs of Equation~\eqref{eq:recdn} have been given by Remmel, Wilf, D\'esarm\'enien and, more recently, Benjamin and Ornstein\footnote{Another bijection is described by Rakotondrajao~\cite[Sec.\ 3]{Rak}, but it appears to be flawed: for $n=5$, both $(5,(13)(24))$ and $(5,(1324))$ seem to be mapped to $(124)(35)$.}.
Remmel's bijection~\cite{Remmel} is quite complicated, and it proves a $q$-analogue of Equation~\eqref{eq:recdn}.
D\'esarm\'enien's bijection~\cite{Desarmenien} first maps derangements to another set of permutations, namely those whose first valley is in an even position.
Wilf's bijection~\cite{Wilf} is easy to program, but it is recursive. Benjamin and Ornstein's bijection~\cite{Benjamin} is perhaps the simplest, but even then, its description requires four different cases.

In this note we present a new, arguably simpler bijective proof of Equation~\eqref{eq:recdn}. 
We describe a bijection 
$\psi:\D^\ast_n\to \F^\ast_n$, where $\D^\ast_n=\D_n\setminus\{(1,2)(3,4)\dots(n-1,n)\}$ and $\F^\ast_n=\F_n$ when $n$ is even, and 
$\D^\ast_n=\D_n$ and $\F^\ast_n=\F_n\setminus\{(1)(2,3)\dots(n-1,n)\}$ when $n$ is odd.

As in~\cite{Benjamin}, we write derangements in cycle notation so that each cycle begins with its smallest element, and cycles are ordered by increasing first element. On the other hand, we write permutations in $\F_n$ with their fixed point at the beginning. 

Let $\pi\in \D^\ast_n$, and let $k$ be the largest non-negative integer such that the cycle notation of $\pi$ starts with $(1,2)(3,4)\dots(2k-1,2k)$. 
Note that $0\le k<n/2$, since $\pi\neq(1,2)(3,4)\dots(n-1,n)$.
To define $\psi(\pi)\in\F^\ast_n$, consider two cases:
\begin{enumerate}[(i)]
\item If the cycle containing $2k+1$ has at least 3 elements, change the first $k+1$ cycles of $\pi$ as follows:
\begin{align*}
\pi&=(1,2)(3,4)\dots(2k-1,2k)(2k+1,a_1,a_2,\dots,a_j)\dots\\
\psi(\pi)&=(1)(2,3)(4,5)\quad\dots\quad(2k,a_1)(2k+1,a_2,\dots,a_j)\dots
\end{align*}
Note that, if $k=0$, then $\{1,2,\dots,2k\}=\emptyset$ and the fixed point in $\psi(\pi)$ is $a_1$.
\item Otherwise, change the first $k+2$ cycles of $\pi$ as follows:
\begin{align*}
\pi&=(1,2)(3,4)\dots(2k-1,2k)(2k+1,a_1)(2k+2,a_2,\dots,a_j)\dots \\
\psi(\pi)&=(1)(2,3)(4,5)\quad\dots\quad(2k,2k+1)(2k+2,a_1,a_2,\dots,a_j)\dots
\end{align*}
\end{enumerate}

The inverse map $\psi^{-1}$ has a similar description. Given $\sigma\in \F^\ast_n$, let $\ell$ be the fixed point of $\sigma$, and consider two cases. If $\ell\neq 1$, merge the cycles containing $\ell$ and $1$ as follows:
$$\sigma=(\ell)(1,a_2,\dots,a_j)\dots\quad \mapsto\quad\psi^{-1}(\sigma)=(1,\ell,a_2,\dots,a_j)\dots$$
Otherwise, let $\sigma'$ be the derangement of $\{2,\dots,n\}$ obtained by removing the fixed point $1$ from~$\sigma$, apply $\psi$ to $\sigma'$ (using the above description, but identifying $\{2,\dots,n\}$ with $\{1,\dots,n-1\}$ in an order-preserving fashion), and replace its fixed point $(\ell)$ with the $2$-cycle $(1,\ell)$ to get $\psi^{-1}(\sigma)$.

As an example, below are the images by $\psi$ of all derangements in $\D_4$ and some in $\D_5$, with the entry $a_1$ colored in \red{red} in case (i) and in \blue{blue} in case (ii):
$$
\begin{array}{c|c|c|c|c|c|c|c|c|c}
\pi & (12)(34) & (1\blue3)(24) & (1\blue4)(23) & (1\red234) & (1\red243) & (1\red324) & (1\red342) & (1\red423) & (1\red432) \\ \hline
\psi(\pi) & -     & (1)(2\blue34) & (1)(2\blue43) & (\red2)(134) & (\red2)(143) & (\red3)(124) & (\red3)(142) & (\red4)(123) & (\red4)(132) 
\end{array}
$$
$$
\begin{array}{c|c|c|c|c|c|c|c|c}
\pi         & -                 & (12)(3\blue45)  & (12)(3\blue54)  & (1\red23)(45) & (1\blue3)(245) & (1\blue4)(235) & (1\red54)(23) &\dots\\ \hline
\psi(\pi) & (1)(23)(45)  & (1)(2\blue4)(35) & (1)(2\blue5)(34) & (\red2)(13)(45)  &(1)(2\blue345) & (1)(2\blue435) & (\red5)(14)(23) & \dots 
\end{array}
$$

\bibliographystyle{plain}
\bibliography{bijectionderangements}

\end{document}